\documentclass[11pt]{amsart}
\usepackage{amssymb}
\numberwithin{equation}{section}

\newtheorem{theorema}{Theorem}
\newtheorem{theorem}{Theorem}[section]

\newtheorem{corollary}[theorem]{Corollary}
\newtheorem{proposition}[theorem]{Proposition}

\theoremstyle{definition}
\newtheorem{definition}[theorem]{Definition}
\theoremstyle{remark}
\newtheorem{remark}[theorem]{Remark}

\newcommand{\R}{\mathbb{R}}
\newcommand{\C}{\mathbb{C}}
\newcommand{\Z}{\mathbb{Z}}
\newcommand{\Q}{\mathbb{Q}}
\newcommand{\bbT}{\mathbb{T}}
\newcommand\lie[1]{\mathfrak{#1}}

\newcommand{\fh}{\lie{h}}
\newcommand{\fg}{\lie{g}}

\newcommand{\fk}{\lie{k}}

\def    \inv    {^{-1}}

\begin{document}

\title{Toric integrable metrics on tori are flat}
\author{Eugene Lerman}
\author{Nadya Shirokova}
\address{Department of
Mathematics, University of Illinois, Urbana, IL 61801}
\email{lerman@math.uiuc.edu}
\email{nadya@math.uiuc.edu}
\date{\today}
\thanks{Partially supported by NSF grant DMS - 980305 and 
the American Institute of Mathematics.}

\begin{abstract} By studying completely integrable torus actions on contact 
manifolds we prove a conjecture of Toth and Zelditch that toric
integrable geodesic flows on tori must have flat metrics.
\end{abstract}

\maketitle

\section{Introduction} In a recent paper Toth and Zelditch 
studied the relation between the dynamics of the geodesic flow on a
compact Riemannian manifold $(Q, g)$ and the growth rate of $L^\infty$
norms of $L^2$-normalized eigenfunctions of the Laplace operator
$\Delta_g$ \cite{TZ}.  They showed that if the square root of the
Laplace operator $\sqrt{\Delta_g}$ is ``quantum completely
integrable'' and has uniformly bounded eigenfunctions then the metric
$g$ is flat, and hence by the Bieberbach theorems $Q$ is finitely
covered by a torus. The proof is particularly transparent when the
geodesic flow is {\bf toric integrable}.  The latter means that there
is an effective action of a torus $\bbT^n$, $n= \dim Q$, on the
punctured cotangent bundle $T^* Q
\smallsetminus Q$ of $Q$ which
\begin{enumerate}
\item commutes with dilations $\rho _t : T^*Q \smallsetminus Q \to 
T^*Q \smallsetminus Q $, $\rho(q, p) = (q, e^t p)$,
\item preserves the standard symplectic form on $T^*Q$ and
\item preserves the energy function $h (q, p) = g^*_q (p, p)$, where $g^*$ 
denotes the metric on $T^*Q$ dual to $g$. (The Hamiltonian flow of $h$
is the geodesic flow.)
\end{enumerate} 
Note that any symplectic group action on the punctured cotangent
bundle which commutes with dilations preserves the Liouville 1-form
and is, therefore, Hamiltonian.  Consequently if a metric on a
manifold $Q$ is toric integrable, the pull-back metric on a finite
cover of $Q$ is toric integrable as well.  One is therefore lead to
wonder if in the case of tori the boundedness of eigenfunctions is
necessary for the flatness of the metric or if toric integrability by
itself is enough.  The main goal of this paper is to prove that, as
conjectured by Toth and Zelditch in \cite{TZ}, toric integrable
metrics on tori are flat:

\begin{theorema}   \label{theoremTZ1}
Suppose that $g$ is a toric integrable metric on a torus $\bbT^n :=
\R^n/\Z^n$.  Then $g$ is flat.
\end{theorema}

The term ``toric integrable'' is apparently due to Toth and Zelditch
(but the objects involved have been studied since the late 70's,
e.g. by Colin de Verdi\`ere, Duistermaat and Guillemin).  It describes a
class of completely integrable systems more general than the systems
with global action-angle variables.  On the other hand toric
integrable systems are much more manageable than arbitrary completely
integrable systems, and one can use the tools of symplectic and
contact geometry to investigate them.

Recall that a symplectic cone is  a symplectic manifold $(M,
\omega)$ with a free proper action $\rho_t$ of the real line which
expands the symplectic form exponentially: ${\rho_t}^* \omega = e^t
\omega$. For example  the punctured cotangent bundle $T^*Q \smallsetminus Q $ 
with the standard symplectic form is a symplectic cone: the real line
acts by dilations $\rho_t (q, p) = (q, e^t p)$ for all $q\in Q$, $p\in
T_q^*Q$.  Given a symplectic cone $(M, \omega , \rho_t)$, a function
$h\in C^\infty (M)$ is {\bf toric integrable} if there is an effective
action of a torus $G$ with $\dim G =
\frac{1}{2} \dim M$ which preserves the symplectic form $\omega$ and
the function $h$ and commutes with dilations $\rho_t$.  Any action of
a torus $G$ on a a symplectic cone $(M, \omega, \rho_t)$ that commutes
with dilations preserves a 1-form $\alpha$ with $d\alpha = \omega$ and
hence is Hamiltonian. Thus toric integrability of a function $h$ on a
symplectic cone $M$ amounts to the existence of $n = \frac{1}{2} \dim
M$ functions $f_1, \ldots, f_n$ which are homogeneous with respect to
dilations, Poisson commute with each other and with $h$ and whose
Hamiltonian flows are all $2\pi$-periodic.

Recall that if $\{f_1, \ldots, f_n\}$ is a completely integrable
system on a symplectic manifold $(M, \omega)$ and if the fibers of the
map $f = (f_1, \ldots f_n) : M \to \R^n$ are compact, then, by
Arnold-Liouville theorem, in a neighborhood of every point of $M$ the
Hamiltonian vector fields of the functions $f_1, \ldots f_n$ generate
a Hamiltonian action of the $n$-torus $\bbT^n$ \cite{Arnold}.
According to Duistermaat there are obstructions to these ``local''
$\bbT^n$ actions to patch up to an action of $\bbT^n$ on $M$ --- the
monodromy of the period lattice \cite{Duis, GSbook}.  Strictly
speaking Duistermaat only considered patching together free torus
actions, but essentially the same argument works in general
\cite{Bouc-Mol}.  If these ``local'' $\bbT^n$ actions patch to a
global $\bbT^n$ action, there is a further obstruction to the
existence of global action-angle variables: the ``Chern class of the
fibration $f: M\to \R^n$.'' (The expression is in quotation marks
because if the torus action is not free then $f$ is not a fibration.
None  the less one can still speak of the ``Chern class'' of $f$
\cite{Bouc-Mol}.)  The second obstruction is easily seen to be
nontrivial --- there are completely integrable systems with global
torus actions but no global action-angle variables.  See, for example,
\cite{Bates}.  Thus toric integrability is weaker than the existence of
global homogeneous action variables.

Toric integrable systems have not been studied systematicly.  We will
see in this paper that a good way to understand them is through the
study contact toric manifolds.  It appears that toric integrability is
rare.  It would be interesting to classify all compact manifolds
admitting toric integrable geodesic flows.  In particular it would be
interesting to find out if there are manifolds other that $S^2$ and
tori that admit such flows.  (Toric integrable metrics on $S^2$ other
than the round one were described by Colin de Verdi\`ere \cite{CdV}.)
This will be addressed elsewhere.  See \cite{ctm} for a first step in
that direction.

\subsection*{Acknowledgments}  We thank Steve Zelditch for suggesting the 
problem and the referee for a number of comments.

\subsection*{A note on notation}  

Throughout the paper the Lie algebra of a Lie group denoted by a
capital Roman letter will be denoted by the same small letter in the
fraktur font: thus $\fg$ denotes the Lie algebra of a Lie group $G$
etc.  The identity element of a Lie group is denoted by 1.  The
natural pairing between $\fg$ and $\fg^*$ will be denoted by $\langle
\cdot, \cdot \rangle$.

When a Lie group $G$ acts on a manifold $M$ we denote the action by an
element $g\in G$ on a point $x\in M$ by $g\cdot x$; $G\cdot x$ denotes
the $G$-orbit of $x$ and so on.  The vector field induced on $M$ by an
element $X$ of the Lie algebra $\fg$ of $G$ is denoted by $X_M$.  The
isotropy group of a point $x\in M$ is denoted by $G_x$; the Lie
algebra of $G_x$ is denoted by $\fg_x$ and is referred to as the
isotropy Lie algebra of $x$.  We recall that $\fg_x = \{ X \in \fg\mid
X_M (x) = 0\}$.

If $P$ is a principal $G$-bundle then $[p, m]$ denotes the point in the
associated bundle $P\times _G M = (P\times M)/G$ which is the orbit of
$(p,m) \in P\times M$.

\section{From toric integrability to contact toric manifolds}

We now start the proof of Theorem~\ref{theoremTZ1}.  As the first
step, following Toth and Zelditch, let us reduce the proof of
Theorem~\ref{theoremTZ1} to a statement about actions of tori on their
punctured cotangent bundles.  To wit, suppose we know that any action
of an $n$-torus $G$ on $M = T^* \bbT^n \smallsetminus \bbT^n$ which is
symplectic and commutes with the action of $\R^+$, is actually a free
action.  Then, as indicated in \cite{TZ} we can apply a theorem of
Ma\~ne
\begin{theorema}[Ma\~ne, \cite{M}]\label{theoremM}
Let $(Q, g)$ be a Riemannian manifold with a geodesic flow $\phi_t :
T^*Q \smallsetminus Q \to T^*Q \smallsetminus Q$.  Suppose the flow
$\phi_t$ preserves the leaves of a non-singular Lagrangian foliation
of $T^*Q \smallsetminus Q$.  Then $(Q, g)$ has no conjugate points.
\end{theorema}

\noindent
to conclude that the toric integrable metric on $\bbT^n$ has no
conjugate points.  Finally, following \cite{TZ} again, and applying

\begin{theorema}[Burago-Ivanov, \cite{BI}]\label{theoremBI}
A metric on a torus $\bbT^n$ with no conjugate points is flat.
\end{theorema}

\noindent
we can conclude that a toric-integrable metric is flat.  To summarize,
in order to prove Theorem~\ref{theoremTZ1} it is enough to show

\begin{proposition}
Suppose that an $n$ torus $G$ acts effectively and symplecticly on the
punctured cotangent bundle $T^* \bbT^n \smallsetminus \bbT^n$, and suppose
that the action commutes with the action of $\R$.  Then the action
of $G$ is free.
\end{proposition}  
Clearly the lift of the left multiplication on $\bbT^n$ to the
cotangent bundle $T^*\bbT^n$ satisfies both the hypotheses of the
proposition and the conclusion.  The crux is to show that an arbitrary
action of $G\simeq \bbT^n$ satisfying the hypotheses has to be free.

Under the hypotheses of the proposition, the action of $G$ descends to
an action on the co-sphere bundle $M:= (T^* \bbT^n \smallsetminus
\bbT^n)/\R$. Moreover this induced action $G$ preserves the natural
contact structure $\xi$ on $M$ (we'll discuss contact structures in
more detail in the next section).  Consequently the proof of
Theorem~\ref{theoremTZ1} reduces to

\begin{theorema}\label{main_thm}
Suppose that an $n$ torus $G$ acts effectively on the co-sphere bundle
$M:= S(T^* \bbT^n)$ of the standard $n$-torus $\bbT^n$ preserving the
standard contact structure on $M$.  Then the action of $G$ is free.
\end{theorema}

Our strategy is to study completely integrable torus actions on
arbitrary (compact connected co-oriented) contact manifolds and to
show that if an action is not free then the underlying manifold cannot
be the product of a torus and a sphere of the appropriate dimensions.

\section{Group actions on contact manifolds}

Recall that a co-oriented contact manifold is a pair $(M, \xi)$ where
$\xi \subset TM$ is a distribution globally given as the kernel of a
1-form $\alpha$ such that $d\alpha |_\xi$ is nondegenerate.  Such a
1-form $\alpha$ is called a contact form and the distribution $\xi$ is
called a contact structure.  A co-sphere bundle $S(T^*N)$ of a
manifold $N$ (defined with respect to some metric) is a natural
example of a contact manifold: the contact form is the restriction of
the Liouville 1-form to the co-sphere bundle.

The condition that a distribution $\xi \subset TM$ is contact is
equivalent to: the punctured line bundle $\xi^\circ \smallsetminus M$
is a symplectic submanifold of the punctured cotangent bundle $T^*M
\smallsetminus M$, where $\xi^\circ $ denotes the annihilator of $\xi$ in
$T^*M$.  Note that if $\xi = \ker \alpha$ then the 1-form $\alpha$ is
a nowhere zero section of the line bundle $\xi ^\circ \to M$.
(Conversely any nowhere vanishing section of $\xi^\circ$ is a contact
form.)  Thus if $\xi =
\ker \alpha$ then $\xi^\circ \smallsetminus M$ has two components.  If a
compact connected Lie group $G$ acts on $M$ and preserves the contact
distribution $\xi$, then the action of $G$ on $\xi^\circ$ maps the
components of $\xi^\circ \smallsetminus M$ into themselves.  Hence
given a contact 1-form $\alpha$ with $\xi = \ker \alpha$, we can
average it over $G$ and obtain a $G$-invariant contact form
$\bar{\alpha}$ with $\xi = \ker \bar{\alpha}$.  Note that each
component of $\xi^\circ \smallsetminus M$ is the symplectization of
$(M, \xi)$.

\begin{definition}
An action of a torus $G$ on a contact manifold $(M, \xi)$ is {\bf
completely integrable} if it is effective, preserves the contact
structure  $\xi$ and if $2\dim G = \dim M +1$.

A {\bf contact toric $G$-manifold} is a co-oriented contact manifold
with a completely integrable action of a torus $G$.
\end{definition}

Note that if an action of a torus $G$ on $(M, \xi)$ is completely
integrable, then the action of $G$ on a component $\xi^\circ_+$ of
$\xi ^\circ \smallsetminus M$ is a completely integrable Hamiltonian
action and thus $\xi_+^\circ$ is a symplectic toric manifold (for more
information on symplectic toric manifolds and orbifolds see \cite{D}
and \cite{LT}).

Completely integrable torus actions on co-oriented contact manifolds
and contact toric manifolds have been studied by Banyaga and Molino
\cite{BM1, BM2, B} and by Boyer and Galicki \cite{BG}.  To state their results
it would be convenient to first digress on the subject of moment maps
for group actions on contact manifolds.

If a Lie group $G$ acts on a manifold $M$ preserving a contact form
$\alpha$, the corresponding {\bf $\alpha$-moment map} $\Psi_\alpha :M
\to \fg^*$  is defined by
\begin{equation}\label{stupid_eq}
\langle \Psi _\alpha (x) , X \rangle = \alpha _x (X_M (x))
\end{equation}
for all $x\in M$ and all $X\in \fg$, where $X_M$ denotes the vector
field corresponding to $X$ induced by the infinitesimal action of the
Lie algebra $\fg$ of the group $G$: $X_M (x) = \frac{d}{dt} |_{t=0}
(\exp tX)\cdot x$.

Note that if $f$ is a $G$-invariant nowhere zero function, then
$\alpha' = f\alpha$ is also a $G$-invariant contact form defining the
same contact structure.  Clearly the corresponding moment map
$\Psi_{\alpha'}$ satisfies $\Psi_{\alpha'} = f\Psi_\alpha$.  Thus the
definition of a contact moment map above is somewhat problematic: it
depends on a choice of an invariant contact form rather then solely on
the contact structure and the action.  Fortunately there is also a
notion of a contact moment map that doesn't have this problem.
Namely, suppose again that a Lie group $G$ acts on a manifold $M$
preserving a contact structure $\xi$ (and its co-orientation).  The
lift of the action of $G$ to the cotangent bundle then preserves a
component $\xi^\circ _+$ of $\xi^\circ \smallsetminus M$.  The
restriction $\Psi = \Phi|_{\xi_+ ^\circ}$ of the moment map $\Phi$ for
the action of $G$ on $T^*M$ to depends only on the action of the group
and on the contact structure.  Moreover, since $\Phi: T^*M \to
\fg^*$ is given by the formula 
$$
\langle \Phi (q, p), X\rangle = \langle p, X_M (q)\rangle 
$$ 
for all $q\in M$, $p\in T^*_q M$ and  $X \in \fg$, we see
that if $\alpha$ is any invariant contact form with $\ker \alpha =
\xi$ then $\langle \alpha^* \Psi (q, p), X\rangle = \langle \alpha ^*
\Phi (q, p), X\rangle = \langle \alpha _q, X_M (q) \rangle = \langle
\Psi _\alpha (q), X\rangle$.  Here we think of $\alpha$ as a section
of $\xi^\circ_+ \to M$.  Thus $\Psi \circ \alpha = \Psi _\alpha$, that
is, $\Psi = \Phi|_{\xi^\circ}$ is a ``universal'' moment map.

There is another reason why the universal moment map $\Psi :\xi^\circ_+
\to \fg^*$ is a more natural notion of the moment map than the one
given by (\ref{stupid_eq}).  The vector fields induced by the action
of $G$ preserving a contact distribution $\xi$ are contact.  The space
of contact vector fields is isomorphic to the space of sections of the
bundle $TM/\xi \to M$.  Thus a contact group action gives rise to a
linear map
\begin{equation}\label{eq**}
\fg \to \Gamma (TM/\xi), \quad X\mapsto X_M \mod \xi .
\end{equation}
The moment map should be the transpose of the map (\ref{eq**}).  The
total space of the bundle $(TM/\xi)^*$ naturally maps into the space
dual to the space of sections $\Gamma (TM/\xi)$: $$ (TM/\xi)^* \ni
\eta \mapsto \left(s \mapsto \langle \eta, s (\pi (\eta))
\rangle\right), $$ where $\pi : (TM/\xi)^* \to M$ is the projection
and $\langle \cdot, \cdot \rangle$ is the paring between the
corresponding fibers of $(TM/\xi)^*$ and $TM/\xi$.  In other words,
the transpose $\Psi: (TM/\xi)^* \to \fg^*$ of (\ref{eq**}) should be
given by
\begin{equation}\label{eq***}
\langle \Psi (\eta), X\rangle = \langle \eta, X_M (\pi (\eta))\!  \mod \xi \rangle
\end{equation}
Under the identification $\xi^\circ \simeq (TM/\xi)^*$, the equation above becomes
$$
\langle \Psi (q, p), X\rangle = \langle p, X_M (q)\rangle
$$ 
for all $q\in M$, $p\in \xi_q^\circ$ and $X\in \fg$, which is
the definition of $\Psi $ given earlier as the restriction to $\xi
^\circ _+$ of the moment map for the lifted action of $G$ on the
cotangent bundle $T^*M$.

Thus part of the above discussion can be summarized as 

\begin{proposition}
Let $(M, \xi)$ be a co-oriented contact manifold with an action of a
Lie group $G$ preserving the contact distribution and its
co-orientation.  Suppose there exists an invariant 1-form $\alpha$
with $\ker \alpha = \xi$.  Then the moment map $\Psi_\alpha $ for the
action of $G$ on $(M, \alpha)$ and the moment map $\Psi$ for the
action of $G$ on the symplectization $\xi^\circ _+$ are related by $$
\Psi \circ \alpha = \Psi _\alpha.
$$ 
Here $\xi^\circ _+$ is the component of $\xi^\circ \smallsetminus
0$ containing the image of $\alpha :M \to \xi ^\circ$.
\end{proposition}

\begin{remark}
We will refer to $\Psi: \xi^\circ_+ \to \fg^*$ as the moment map for
the action of a Lie group $G$ on a co-oriented contact manifold $(M,
\xi = \ker \alpha)$.
It is easy to show that $\Psi$ is $G$-equivariant with respect to the
given action of $G$ on $M$ and the coadjoint action of $G$ on $\fg^*$.
Hence for any invariant contact form $\alpha$ the corresponding moment
map $\Psi_\alpha :M \to \fg^*$ is also $G$-equivariant.
\end{remark}

\begin{definition} 
Let $(M, \xi = \ker \alpha )$ be a co-oriented contact manifold with
an action of a Lie group $G$ preserving the contact distribution and
its co-orientation.  Let $\Psi : \xi ^\circ_+ \to \fg^*$ denote the
corresponding moment map.  We define the {\bf moment cone} $C(\Psi)$
to be the image of a connected component $\xi_+^\circ$ of $\xi^\circ
\smallsetminus M$ plus the origin: 
$$ 
C(\Psi) := \Psi (\xi_+^\circ)\cup \{0\}.  
$$

Remark that 
$$ 
C(\Psi) = \R^+ \Psi_\alpha (M) \cup \{0\}
$$
where $\Psi_\alpha :M \to \fg^*$ is the $\alpha$-moment map.
\end{definition}

Note that the moment cone does not depend on the choice of a contact
form; it is a true invariant of the co-oriented contact structure and
the group action.

\begin{remark}
An action of a Lie group $G$ on a manifold $M$ preserving a contact
form $\alpha$ is completely encoded in the moment map $\Psi_\alpha: M
\to \fg^*$.  Therefore it will be convenient for us to think of a
{\bf contact toric $G$-manifold} as an equivalence class of triples $(M,
\alpha, \Psi_\alpha :M \to \fg^*)$ where the $\Psi_\alpha$ is the
moment map for a completely integrable  action of a torus $G$ on a contact
manifold $(M, \alpha)$, or, somewhat more sloppily, as a triple   $(M,
\alpha, \Psi_\alpha :M \to \fg^*)$.
\end{remark}

\section{Proof of Theorem~\ref{main_thm}}

Banyaga and Molino made the first step towards classifying compact
connected contact toric manifolds in \cite{BM1}.  A revised version of
this paper circulated as the preprint \cite{BM2}.  The main
classification result of \cite{BM2} is cited in \cite{B} roughly as
follows:

\begin{theorem}\label{BMthm}
Let $(M, \alpha, \Psi_\alpha:M \to \fg^*)$ be a compact connected
contact toric $G$-manifold.

Suppose the action of $G$ on $M$ is free.  Then the orbit space $M/G$
is diffeomorphic to a sphere.  If additionally $\dim G > 2$ then the
map $\bar{\Psi}_\alpha : M/G \to \fg^*$ induced by the moment map
$\Psi_\alpha$ is an embedding.  If furthermore $\dim G > 3$, then $M$
is the co-sphere bundle of $G$, i.e., $M = S(T^*G)$.

Suppose the action of $G$ on $M$ is not free and suppose $\dim G> 2$.
Then the moment cone $C(\Psi)$ is a convex polyhedral cone and
the map $\bar{\Psi}_\alpha: M/G \to \fg^*$ induced by the moment map
is an embedding.  Moreover the cone  $C(\Psi)$ determines the
contact toric manifold.
\end{theorem} 
\begin{remark}
It is easy to construct examples of a completely integrable action of
a 2-torus on a contact 3-torus for which the fibers of the
corresponding moment map are not connected: let $M=
\bbT^3$ with coordinates $\theta_1, \theta_2$ and $t$, let 
$\alpha = \cos 2t\, d\theta_1 +\sin 2t\, d\theta_2$ be the contact
form, and let $\bbT^2$ act by $(\mu, \nu) \cdot (\theta_1,
\theta_2, t) = (\theta_1 + \mu, \theta_2 + \nu, t)$. 
Also there are examples of completely integrable 2-torus actions on
overtwisted lens spaces for which the corresponding moment cones are
not convex.  See \cite{L-IJM}.
\end{remark}
Contact toric manifolds have also been studied by Boyer and Galicki
\cite{BG}.  The following result is implicit in their paper:
\begin{theorem}\label{BGthm}
Let $(M, \alpha, \Psi_\alpha: M\to \fg^*)$ be a compact connected
contact toric $G$-manifold.  Suppose there exits a vector $X\in \fg$
such that the component of the moment map $\langle \Psi _\alpha,
X\rangle$ is strictly positive on $M$.  Then $M$ is a Seifert bundle
over a (compact) symplectic toric orbifold.
\end{theorem} 

We remind the reader that a symplectic toric orbifold is a symplectic
orbifold with a completely integrable Hamiltonian torus action.  
Compact connected symplectic toric orbifolds were classified in \cite{LT}.

\begin{proof}[Proof of Theorem~\ref{BGthm}]
Since $M$ is compact, the image $\Psi_\alpha (M)$ is compact.
Therefore the set of vectors $X' \in \fg$, such that the function
$\langle \Psi_\alpha, X'\rangle$ is strictly positive on $M$, is open.
Hence we may assume that $X$ lies in the integral lattice $\Z_G :=
\ker (\exp: \fg \to G)$ of the torus $G$.  Let $H = \{\exp tX \mid
t\in \R\}$ be the corresponding circle subgroup of $G$.

Let $f(x) = 1/(\langle \Psi_\alpha (x) ,X\rangle)$ and let $\alpha' =
f \alpha$.  The form $\alpha '$ is another $G$-invariant contact form
with $\ker \alpha' = \xi$.  The moment map $\Psi_{\alpha'}$ defined by
$\alpha'$ satisfies $\Psi_{\alpha'} = f \Psi_\alpha$.  Therefore
$\langle \Psi_{\alpha'} (x), X\rangle = 1$ for all $x\in M$.
 
Since the function $\langle \Psi_\alpha , X \rangle $ is nowhere zero,
the action of $H$ on $M$ is locally free.  Consequently the induced
action of $H$ on the symplectization $(N, \omega) = (M\times \R,
d(e^t\alpha'))$ is locally free as well.  Hence any $a\in \R $ is a
regular value of the $X$-component $\langle \Phi, X\rangle$ of the
moment map $\Phi$ for the action of $G$ on the symplectization $(N,
\omega)$.  Note that $\Phi (x, t) = -e^t \Psi_{\alpha'}(x)$.  Now
$M\times \{0\}$ is the $-1$ level set of $\langle \Phi, X\rangle$.
Therefore $B:= (\langle \Phi, X\rangle )\inv (-1)/H \simeq M/H$ is a
(compact connected) symplectic orbifold with an effective Hamiltonian
action of $G/H$.  The orbit map $\pi :M\simeq (\langle \Phi, X\rangle
)\inv (-1) \to B$ makes $M$ into a Seifert bundle over $B$.  A
dimension count shows that the action of $G/H$ on $B$ is completely
integrable.
\end{proof}

\begin{remark}
It is easy to see that the moment cone for the action of $G$ on $M$ is
the cone on the moment polytope of $B$.  In particular it is a proper
polyhedral cone, that is, it contains no linear subspaces.
\end{remark}
\begin{corollary}
Let $(M, \alpha, \Psi_\alpha: M \to \fg^*)$ be a (compact connected)
contact toric manifold.  Suppose there exits a vector $X\in \fg$ such
that the component of the moment map $\langle \Psi_\alpha, X\rangle$
is strictly positive on $M$.  Then $\dim _\Q H^1 (M, \Q) \leq 1$.
\end{corollary}

\begin{proof}
By Theorem~\ref{BGthm} the manifold $M$ is a Seifert bundle over a
compact connected symplectic toric orbifold $B$.  A generic component
of the moment map on $B$ is a Morse function with all indices even.
The  Morse inequalities hold rationally for Morse functions on
orbifolds (see \cite{LT}).  Therefore the first cohomology $H^1 (B,
\Q)$ is zero.

Next we apply the Gysin sequence to the map $\pi: M \to B$.  Since the
Gysin sequence comes from the collapse of the Leray-Serre spectral
sequence for $\pi$ and since rationally the "fibration" $\pi$ is a
circle bundle, the Gysin sequence does exist.  We have $0 \to H^1 (B,
\Q) \to H^1 (M, \Q) \stackrel{\pi _*}{\to} H^0 (B, \Q) \to H^2 (B, \Q)
\to \cdots$.  Since $H^0 (B, \Q) = \Q$ and since $H^1 (B, \Q) = 0$, the
result follows. 
\end{proof}
We conclude immediately 

\begin{corollary}\label{easy-cor}
If $M$ is a contact toric manifold satisfying the hypotheses of the
Theorem~\ref{BGthm}, then $M$ is not the co-sphere bundle of a torus.
\end{corollary}

Combining Corollary~\ref{easy-cor} with Theorem~\ref{BMthm} we see
that if an $n$-dimensional torus $G$ ($n>2$) acts on the co-sphere
bundle $M = S(T^*\bbT^n)$ preserving a contact form $\alpha$ and if
the action is not free, then the corresponding moment cone $C(\Psi)$
contains a linear subspace $P$ of positive dimension.

\begin{proposition} Suppose $(M, \alpha , \Psi_\alpha :M \to \fg^*)$ is a 
compact connected contact toric $G$-manifold of dimension $2n-1>3$,
the action of $G$ on $M$ is not free and the moment cone $C(\Psi)$
contains a linear subspace  $P$ of dimension $k$, $0<k<n$. Then $\dim
H^1(M, \Q) = k \not = n =\dim H^1(S(T^* \bbT^n), \Q) $.
\end{proposition}

\begin{proof} By Theorem~\ref{BMthm} the fibers of
the contact moment map $\Psi_\alpha$ are connected.  Let $\Phi:
M\times \R \to \fg^*$ denote the symplectic moment map for the
Hamiltonian action of $G$ on the symplectization $(M\times \R, d(e^t
\alpha))$ of $(M, \alpha)$.  It is given by  $\Phi (m, t)
= -e^t \Psi _\alpha (m)$.  Thus $\Phi (M\times \R)\cup \{0\} = -C(\Psi
)$ and the fibers of $\Phi$ are connected.  The triple $(M\times \R,
d(e^t \alpha), \Phi )$ is a symplectic toric manifold.  Since the
image of $\Phi$ is contractible and since the fibers of $\Phi$ are
connected, it follows from a result of Lerman, Tolman and Woodward
(Lemma~7.2 and Proposition~7.3 in \cite{LT}), that the image of $\Phi$
determines the symplectic toric manifold $(M\times \R, d(e^t \alpha),
\Phi )$ uniquely.\footnote{In proving the result Lerman, Tolman and
Woodward rediscovered the ideas of Boucetta and Molino
\cite{Bouc-Mol}} In particular the image determines the first
cohomology group of $M$.

A standard argument for Hamiltonian $G$-spaces implies that the
subspace $P$ is the annihilator of the Lie algebra $\fh$ of a subtorus
$H$ of $G$.  Since $H$ is a subtorus, the exact sequence 
$$
1\to H \to G\to G/H \to 1$$ 
splits.  Let $K\simeq G/H$ be a subtorus in $G$
complementary to $H$ and let $\fg^* = \fh^* \times \fk^*$ be the
corresponding splitting of the duals of the Lie algebras.  Then
$P\simeq \fk^*$ and $C(\Psi_\alpha) = P\times W$ where $W\subset
\fh^*$ is a proper cone.  It follows from a theorem of Delzant \cite{D} 
that the exists a basis $w_1, w_2,\ldots, w_l$ of weight lattice of
the torus $H$ so that the edges of the cone $W$ are of the form $\R^+
w_i$.  In particular the representation of $H$ on $\C^l$ defined by
the infinitesimal characters $w_1, \ldots, w_l$ has the property that
the image of the corresponding moment map is $W$.  Consequently we can
realize $\Phi (M\times \R)$ as the image of the moment map for
the product action of $K\times H$ on $K\times \left(\fk^* \times \C^l
\smallsetminus (0,0)\right)\subset T^*K \times \C^l$.  Therefore $M\times
\R$ is $G$-equivariantly symplectomorphic to $K\times \left(\fk^* \times \C^l
\smallsetminus (0,0)\right)$, which is homotopy equivalent to 
$K\times S^{k+ 2l -1}= (S^1)^k\times S^{k+ 2l -1} $.   Since
$0< k < n$ and since $l = n-k >0$, $k = \dim H^1 ((S^1)^k\times S^{k+ 2l
-1}) = \dim H^1 (M)$.
\end{proof}

It remains to consider the case where $(M, \alpha, \Psi_\alpha :M\to
\fg^*)$ is a contact toric manifold of dimension $3$ and the action of
the 2-torus $G$ is not free. 

\begin{proposition}\label{prop_lens}  Suppose $(M, \alpha, \Psi_\alpha :M\to
\fg^*)$ is a compact connected contact toric $G$-manifold of dimension $3$ and 
suppose the action of
the 2-torus $G$ is not free.  Then $M$ is a lens space and hence
cannot be the co-sphere bundle of a 2-torus.
\end{proposition}
\begin{remark}
We consider the 3-manifold $S^1 \times S^2$ a lens space. 
\end{remark}

\begin{proof}[Proof of Proposition~\ref{prop_lens}]
As we remarked earlier the symplectization $(M\times \R, d(e^t
\alpha), \Phi (m,t) = - e^t \Psi_\alpha (m))$ is a symplectic toric
manifold.  Delzant \cite{D} showed that for symplectic toric manifolds
all the isotropy groups are connected and all fixed points are
isolated.  If a point $x\in M$ is fixed by the action of $G$ then the
line $\{x\}\times \R$ is fixed by the action of $G$ on $M\times \R$.
Therefore the action of $G$ on a contact toric $G$-manifold has no
fixed points (one can also give a direct proof of this fact).

Next we use the fact that $\dim M = 3$ and $\dim G = 2$.  By the above
observations the isotropy groups for the action of $G$ on $M$ are
either trivial or circles.  If the isotropy group of $x\in M$ is
trivial, then a neighborhood of the orbit $G\cdot x$ in $M$ is
$G$-equivariantly diffeomorphic to $G\times (-\epsilon, \epsilon)$ for
some small epsilon.  If the isotropy group of a point $x\in M$ is a
circle $H < G$, then a neighborhood of the orbit $G\cdot x$ in $M$ is
$G$-equivariantly diffeomorphic to $G\times _H D^2$ for a small disk
$D^2 = \{z \in \C \mid |z| < \epsilon\}$.  Moreover the action of $H$
on $D^2$ must be effective; hence we may identify $H$ with $S^1$ in such a way 
that the action of $H$ on $D^2$ is given by $\lambda \cdot z = \lambda z$.

We conclude that  locally  $M/G$ is homeomorphic to either $(-\epsilon, \epsilon)$ or to  $D^2 /S^1 \simeq [0,\epsilon)$.  Thus if the action of $G$ on $M$ is not free, then $M/G$ is
a one-dimensional manifold with boundary.  Since  $M$ is compact and 
connected, $M/G$ has to be an interval.  Therefore, by a theorem of
Haefliger and Salem (Proposition~4.2 in \cite{HS}), $M$ as a $G$-space
is uniquely determined by the isotropy representations at the points
in $M$ above the endpoints of the interval $M/G$.  It is easy to see
that in this case $M$ is diffeomorphic to two solid tori glued
 along their  boundaries.  We
conclude that if $M$ is a three dimensional compact connected contact
toric manifold and if the action of a 2-torus $G$ is not free, then
$M$ is a lens space.  In particular $M$ is not diffeomorphic to
$\bbT^3 = S(T^* \bbT^2)$.
\end{proof}

This finishes the proof of the main result, Theorem~\ref{main_thm} and
consequently the proof of Theorem~\ref{theoremTZ1}.

\end{document}